\documentclass[]{article}
\usepackage{graphicx,amsfonts,amsmath,amsthm,amssymb,amscd,epic,eepic,rotating}

\def\tr{{\rm tr}}
\def\PSL{{\rm PSL}(2,{\mathbb C})}

\def\U{{\cal U}}
\def\RP{{\cal RP}}
\def\D{{\cal D}}

\newtheorem{theorem}{Theorem}[section]

\newtheorem{proposition}{Proposition}[section]
\newtheorem{cor}{Corollary}[section]
\newtheorem{remark}{Remark}[section]

\title{Two-generator Kleinian orbifolds}
\author
{Elena Klimenko
\thanks{Gettysburg College, Mathematics Department, 300 N Washington St.,
CB 402, Gettysburg, PA 17325, USA; e-mail
{\tt yklimenk@gettysburg.edu}}
\thanks{supported by Gettysburg College
Research and Professional Development Grant, 2005--2006}
\and
Natalia Kopteva
\thanks{ LATP, UMR CNRS 6632, CMI,  Universit\'e de Provence,
39 rue F. Joliot Curie, 13453 Marseille cedex 13, France; e-mail
{\tt kopteva@cmi.univ-mrs.fr}}
\thanks{supported by FP6 Marie Curie IIF Fellowship}
}

\date{\today}

\begin{document}

\maketitle

\begin{abstract}
We give a complete list of orbifolds uniformised by 
discrete non-elementary
two-generator subgroups of $\PSL$ without invariant plane
whose generators and their commutator have real traces.
\end{abstract}

\footnotesize
\noindent
{\bf Mathematics Subject Classification (2000): }
Primary: 57M12;

\noindent
Secondary: 30F40, 22E40, 57M50.

\noindent
{\bf Key words: }
hyperbolic orbifold, Kleinian group, discrete group,
hyperbolic geometry.
\normalsize

\section{Introduction}

Let $\Gamma$ be a non-elementary Kleinian group, and let
$\Omega(\Gamma)$ be the discontinuity set of $\Gamma$.
Following \cite{BP01}, we say that the {\it Kleinian orbifold}
$Q(\Gamma)=({\mathbb H}^3\cup\Omega(\Gamma))/\Gamma$
is an orientable $3$-orbifold with
a complete hyperbolic structure
on its interior ${\mathbb H}^3/\Gamma$ and a conformal structure
on its boundary $\Omega(\Gamma)/\Gamma$.

First, we define the class of $\cal RP$ {\it groups} (two-generator
groups with real parameters) by
$$
{\cal RP}=\lbrace\Gamma : \Gamma=\langle f,g\rangle
{\rm \ for\ some\ } f,g\in{\rm PSL}(2,{\mathbb C})
{\rm \ with\ }
\beta,\beta',\gamma\in {\mathbb R}\rbrace,
$$
where $\beta=\beta(f)={\rm tr}^2f-4$, $\beta'=\beta(g)={\rm tr}^2g-4$,
$\gamma=\gamma(f,g)={\rm tr}[f,g]-2$.
A pair $(f,g)$ such that $(\beta,\beta',\gamma)\in{\mathbb R}^3$
is called an {\it $\RP$ pair}.

Now denote by $\cal D$ the class of non-elementary Kleinian $\RP$ groups
without invariant plane whose generators have real traces.
The purpose of this paper is to describe the Kleinian
orbifolds for the class~$\cal D$.

Before we give an equivalent definition of
 $\cal D$ in terms of parameters,
we recall that an element $f\in\PSL$ with real $\beta=\beta(f)$
is
{\it elliptic}, {\it parabolic}, {\it hyperbolic}, or
{\it $\pi$-loxodromic}
according to whether
$\beta\in[-4,0)$, $\beta=0$, $\beta\in(0,+\infty)$, or
$\beta\in(-\infty,-4)$.
If $\beta\notin[-4,\infty)$, i.e. ${\rm tr} f$ is not real,
then $f$ is called {\it strictly loxodromic}.
Among all
strictly loxodromic elements, only $\pi$-loxodromics have
real $\beta$.

The parameter $\gamma$ is responsible for the mutual position
of invariant planes of $f$ and $g$. For the geometric meaning of
$\gamma$ we refer the reader to~\cite{KK02} and~\cite{KK04-1}.
In particular, from~\cite[Theorem~4]{KK02} we easily deduce the
following

\begin{proposition}\label{spatial_groups}
Let $\Gamma$ be a two-generator subgroup of $\PSL$ and let $f$ and $g$
be its generators with parameters $(\beta,\beta',\gamma)$.
Then $\beta>-4$, $\beta'>-4$,
$\gamma<-\beta\beta'/4$, and $\gamma\not=0$
if and only if all of the following hold:
\begin{enumerate}
\item $\Gamma$ is an ${\cal RP}$ group;
\item $\Gamma$ is non-elementary;
\item $\Gamma$ has no invariant plane (in particular, $\Gamma$ is not
Fuchsian);
\item each of the generators $f$ and $g$ is
either elliptic, parabolic, or
hyperbolic.
\end{enumerate}
\end{proposition}

If we denote by $\cal K$ the class of all Kleinian groups, i.e. discrete
subgroups of $\PSL$, then Proposition~\ref{spatial_groups} gives us 
the description of the class~$\cal D$ in terms of parameters:
{\setlength\arraycolsep{2pt}
\begin{eqnarray*}
{\cal D}=\{\Gamma\in{\cal K} &:& \Gamma=\langle f,g\rangle
{\rm \ for\ some\ } f,g\in\PSL \\
&&{\rm \ with\ }
\beta>-4,\beta'>-4,\
\gamma<-\beta\beta'/4, {\rm \ and\ } \gamma\not=0\}.
\end{eqnarray*}}
The class $\cal D$ was studied in \cite{Kli89}--\cite{KK05}.
To lighten the reading,
we collect the precise references to our previous results
in Table~\ref{refs}.

\begin{table}[htbp]
\begin{center}
\caption{References for the class ${\cal D}$}\label{refs}
{\tabcolsep=1mm
\begin{tabular}{|c|c||c|c|}
\hline
$\beta$  and $f$ & $\beta'$ and $g$ & discreteness & group presentations \\
 & & criteria & and orbifolds\\
\hline\hline
$(-4,0)$ & $(-4,0)$ & Appendix: Theorem~\ref{two_ell} & Appendix: \\
{\footnotesize elliptic} & {\footnotesize elliptic} &
         and Corollary~\ref{cor1}; & Corollary~\ref{cor2}\\
 & & \cite[Theorem]{Kli89}& \\
\hline
 $0$ & $(-4,+\infty)$ &  & \\
 {\footnotesize parabolic} & {\footnotesize elliptic, parabo-} &
  \cite[Theorem 2.3]{KK04-1} & \cite[Theorem~3.1]{KK04-1} \\
    & {\footnotesize lic, hyperbolic} & & \\
\hline
$(0,+\infty)$ & $(0,+\infty)$ & \cite[Theorem 2.1]{Kli01} & 
\cite[Section~3]{Kli01}\\
{\footnotesize hyperbolic} & {\footnotesize hyperbolic} & & \\
\hline
$(-4,0)$ & $(0,+\infty)$ & \cite[Theorem 1]{Kli90}, & [14,
Proposition~3.1,\\
{\footnotesize elliptic} & {\footnotesize hyperbolic} &
\cite[Theorems A and B]{KK02}, &  Proposition 3.2, and\\
 & &\cite[Theorems A and B]{KK04} & Proposition 3.5]\\
\hline
\end{tabular}}
\end{center}
\end{table}

An elliptic element $f$ of order $n$ is said to be {\it non-primitive}
if it is a rotation through $2\pi q/n$, where $q$ and $n$ are
coprime and $1<q<n/2$. If $f$ is a rotation through $2\pi/n$, then it is
called {\it primitive}.

In this paper, we call a pair $(f,g)$ of generators {\it primitive} if
each of $f$ and $g$ is either hyperbolic, or parabolic, or
a primitive elliptic element. It is clear that every $\Gamma\in\D$
has at least one primitive $\RP$ pair of generators.
\medskip

The paper is organised as follows. 
Theorem~\ref{allgroups} in Section~2 gives a presentation
for each $\Gamma\in\cal D$.
Theorem~\ref{th_table} lists the
parameters $(\beta,\beta',\gamma)$ for all primitive $\RP$ pairs
of generators for each $\Gamma\in\D$. The list of parameters 
(without group presentations) first
appeared in \cite{KK04}.

Theorem~\ref{orb_th} in Section~3 lists all orbifolds $Q(\Gamma)$
with $\Gamma\in\cal D$. At the end of Section~3 we give a complete 
list of finite covolume groups in class $\D$
(cocompact and not).

In Appendix, we concentrate on
primitive $\RP$ pairs of elliptics.
We give necessary and sufficient conditions for discreteness of
the non-elementary groups without invariant plane 
generated by such $\RP$ pairs
(Theorem~\ref{two_ell}), reformulate this theorem in terms of parameters
(Corollary~\ref{cor1}),
and describe corresponding orbifolds and group presentations
(Corollary~\ref{cor2}).

Theorem~\ref{two_ell} was proved in \cite{Kli89}, but there is
a mathematical disorder in the English translation of
the original result. So, we reproduce the proof here and complete
the case of two elliptic generators by giving parameters for
such groups and describing corresponding orbifolds.
On the other hand, Appendix can be regarded as an example that
helps to follow the general line of the study of $\RP$ groups.

We remark that the case of two elliptic generators whose
commutator is also elliptic was studied in~\cite{MM03}.

\section{Combinatorics of the class ${\cal D}$: group presentations}

In this section we will give a presentation for every 
group $\Gamma\in\cal D$ (Theorem~\ref{allgroups})
and a complete list of parameters $(\beta,\beta',\gamma)$
for primitive $\RP$ pairs of generators for all
$\Gamma\in\D$ (Theorem~\ref{th_table}). 

The parameters for $\RP$ pairs with non-primitive elliptics
can be calculated using parameters for
primitive elliptics as follows.
Suppose that $f$ is non-primitive elliptic of finite order $n$,
i.e., $\beta(f)=-4\sin^2(q\pi/n)$, where $(q,n)=1$, $1<q<n/2$.
Then there exists an integer $r$ so that $f^r$
is primitive of the same order. 
Obviously, $\langle f,g\rangle=\langle f^r,g\rangle$
and $\beta(f^r)=-4\sin^2(\pi/n)$. By \cite{GM94-1},
$\gamma(f^r,g)=(\beta(f^r)/\beta(f))\gamma(f,g)$.

We need the following group presentations:

\begin{enumerate}
\item
$GT[n,m;q]=\langle f,g\,|\,f^n,g^m,[f,g]^q\rangle$
\item
$PH[n,m,q]=
\langle x,y,z\,|\,x^n,y^2,z^2,(xz)^2,[x,y]^m,(yxyz)^q\rangle$
\item
$H[p;n,m;q]=\langle
x,y,s\,|\,s^2,x^n,y^m,(xy^{-1})^p,(sxsy^{-1})^q,(sx^{-1}y)^2\rangle$
\item
$P[n,m,q]=\langle w,x,y,z\,|\,w^n,x^2,y^2,z^2,(wx)^2,(wy)^2,(yz)^2,
(zx)^q,(zw)^m\rangle$
\item
$Tet[p_1,p_2,p_3;q_1,q_2,q_3]=\langle x,y,z\,|\,
x^{p_1},y^{p_2},z^{p_3},
(xy^{-1})^{q_3},(yz^{-1})^{q_1},(zx^{-1})^{q_2}\rangle$.
The group $Tet[2,2,n;2,q,m]$ is
denoted by $Tet[n,m;q]$ for simplicity.
\item
$GTet_1[n,m,q]=\langle
x,y,z\,|\,x^n,y^2,(xy)^m,[y,z]^q,[x,z]\rangle$
\item
$GTet_2[n,m,q]=\langle
  x,y,z\,|\,x^n,y^2,(xy)^m,(xz^{-1}y^{-1}zy)^q,[x,z]\rangle$
\item
${\cal S}_2[n,m,q]=\langle
x,L\,|\,x^n,(xLxL^{-1})^m,(xL^2x^{-1}L^{-2})^q\rangle$
\item
${\cal S}_3[n,m,q]=\langle
x,L\,|\,x^n,(xLxL^{-1})^m,(xLxLxL^{-2})^q\rangle$
\item
$R[n,m;q]=\langle u,v\,|\,(uv)^n,(uv^{-1})^m,[u,v]^q\rangle$
\end{enumerate}
In the presentations 1--10, the exponents $n,m,q,\dots$ may be
integers (greater than~1), $\infty$, or~$\overline\infty$.
We employ the symbols
$\infty$ and $\overline\infty$ in the following way.
If we have relations of the form $w^n=1$, where $n=\overline\infty$,
we remove them
from the presentation (in fact, this means that the
element $w$ is hyperbolic in the Kleinian group). Further, if we keep
the relations $w^\infty=1$, we get a Kleinian group
presentation where parabolics are indicated.
To get an abstract group presentation, we need to remove such
relations as well.

We assume that $\overline\infty>\infty>x$ and 
$x/\infty=x/\overline\infty=0$
for every real $x$;
$\infty/x=\infty$ and $\overline\infty/x=\overline\infty$
for every positive real $x$;
in particular, $(\infty,n)=(\overline\infty,n)=n$ for every
positive integer $n$.

\begin{theorem}\label{allgroups}
$\Gamma\in{\cal D}$ if and only if $\Gamma$ is isomorphic to one of the following
groups:
\begin{itemize}
\item[1.]
$GT[n,m;q]$,
where $3\leq n\leq m$, and $\cos(\pi/(2q))>\sin(\pi/n)\sin(\pi/m)$;
\item[2.]
$PH[n,m,q]$,
where $4\leq n<\infty$, $n$ is even, $2/n+1/m<1$, and $q\geq 3$ is odd;
\item[3.]
$H[p;n,m;q]$,
where $[p;n,m;q]$ is one of the following:

$[2;2,3;5]$; $[2;2,5;3]$; $[2;3,3;p]$, where $p\geq 5$ is odd;

$[2;3,n;2]$, where $n\geq 5$, $(n,6)=1$;
\item[4.]
$P[n,m,q]$,
where $4\leq n< \infty$, $n$ is even,
$1/n+1/m<1/2$,  $m,q\geq 3$ are odd;
\item[5.]
$Tet[2,3,3;2,3,m]$,
where $m\geq 4$, $(m,3)=1$; or

$Tet[n,m;q]$, where $3\leq n\leq m$, $q\geq 3$ is odd and
$\cos(\pi/q)>\sin(\pi/n)\sin(\pi/m)$;
\item[6.]
$GTet_1[n,m,q]$, where 

either $4\leq n<\infty$, $n$ is even,
$m$ is odd, $1/n+1/m<1/2$,  $q\geq 2$, 

or $n\geq 7$ is odd, $m=3$, $q=2$;
\item[7.]
$GTet_2[n,m,q]$,
where $n,m\geq 3$ are odd, $1/n+1/m<1/2$, and $q\geq 2$;
\item[8.]
${\cal S}_2[n,m,q]$,
where $4\leq n< \infty$, $n$ is even, $2/n+1/m<1$, and $q\geq 2$;
\item[9.]
${\cal S}_3[n,m,q]$,
where $n\geq 3$ is odd, $2/n+1/m<1$, and $q\geq 2$;
\item[10.]
$R[n,2;2]$,
where $n\geq 5$ and $(n,6)=1$.
\end{itemize}
\end{theorem}

{\it Proof} easily follows from Theorem~\ref{th_table}.\qed

\begin{remark}
{\rm
Not all groups in Theorem~\ref{allgroups} are distinct. For example,
both $GT[n,\infty;\overline\infty]$ and ${\cal S}_2[n,\infty,\infty]$
are isomorphic to ${\mathbb Z}_n*{\mathbb Z}$.
}
\end{remark}

If $f\in\PSL$ is a loxodromic element with translation length $d_f$ and
rotation angle $\theta_f$, then
$$
\tr^2f=4\cosh^2 \frac{d_f+i\theta_f}2
$$
and $\lambda_f=d_f+i\theta_f$ is called the
{\it complex translation length} of~$f$.

Note that if $f$ is hyperbolic then $\theta_f=0$ and $\tr^2f=4\cosh^2(d_f/2)$.
If $f$ is elliptic then $d_f=0$ and $\tr^2f=4\cos^2(\theta_f/2)$.
If $f$ is parabolic then $\tr^2f=4$;
by convention we set $d_f=\theta_f=0$.

We define the set
$$
\mathcal{U}=\{u:u=i\pi/p
{\rm \ for\ some\ } p\in{\mathbb Z}, p\geq 2\}\cup[0,+\infty).
$$
In other words, the set $\U$ consists of
all complex translation half-lengths
$u=\lambda_f/2$ for hyperbolic,
parabolic, and primitive elliptic elements~$f$.
Furthermore, we define a function
$t:\mathcal{U}\to\{2,3,4,\dots\}\cup\{\infty,\overline\infty\}$
as follows:
$$
t(u)=\left\{
\begin{array}{lll}
p & {\rm if} & u=i\pi/p,\\
\infty & {\rm if} & u=0,\\
\overline\infty & {\rm if} & u\in(0,+\infty).
\end{array}
\right.
$$

Given $u\in\U$ and $f$ with $\tr^2f=4\cosh^2u$,
$t(u)$ determines the type of $f$ and, moreover, its
order if $f$ is elliptic.
Note also that since we regard $\infty/n=\infty$ and
$\overline\infty/n=\overline\infty$, an expression of the form
$(t(u),n)=1$ with $n>1$ means, in particular, that $t(u)$ is finite
(see Remark~\ref{rem_funcmu} in Appendix as an example).

Now we are ready to introduce Table~\ref{tabgroups}.

\begin{theorem}\label{th_table}
Table~\ref{tabgroups} gives a complete list of parameters 
$(\beta,\beta', \gamma)$ for primitive $\RP$ pairs $(f,g)$ of
generators for all $\Gamma\in\D$. Moreover, for each triple 
$(\beta,\beta',\gamma)$ a presentation of the group $\Gamma$
is given.
\end{theorem}

\begin{proof}
In order to prove the theorem, we need only to summarise the results of our
previous papers (see Table~\ref{refs}).
\end{proof}

\begin{sidewaystable}
\centering
\caption{The parameters and presentations for
the groups $\Gamma=\langle f,g\rangle\in\D$ so that
$(f,g)$ is a primitive $\RP$ pair}\label{tabgroups}
{\footnotesize
\begin{tabular}{rcccc}
\multicolumn{5}{c}{$u,v,w\in\U$ and $n,m,p,q$ are positive integers}\\
\hline
\rule[-2ex]{0ex}{5ex}  & $\beta=\beta(f)$ & $\gamma=\gamma(f,g)$ & $\beta'=\beta(g)$ & $\Gamma=\langle f,g\rangle$ \\
\hline
1 &  \rule[-2ex]{0ex}{6ex} $4\sinh^2u$, $t(u)\geq 3$
  & $-4\cosh^2w$, $(t(w),2)=2$,
  & $4\sinh^2v,$ $t(v)\geq 3$
  & $GT[t(u),t(v);t(w)/2]$ \\
  & \multicolumn{3}{c}{$\cos\frac{\pi}{t(w)}>
     \sin\frac{\pi}{t(u)}\sin\frac{\pi}{t(v)}$}& \\
2 &  \rule[-2ex]{0ex}{6ex} $4\sinh^2u$, $t(u)\geq 3$
  & $-4\cosh^2w$, $(t(w),2)=1$,
  & $4\sinh^2v,$ $t(v)\geq 3$
  & $Tet[t(u),t(v);t(w)]$ \\
  & \multicolumn{3}{c}{$\cos\frac{\pi}{t(w)}>
     \sin\frac{\pi}{t(u)}\sin\frac{\pi}{t(v)}$}& \\
3 & \rule[-2ex]{0ex}{6ex} $-4\sin^2\frac{\pi}n$, $n\geq 5$, $(n,2)=1$
  & $-(\beta+2)^2$
  & $4(\beta+4)\cosh^2u-4$, $t(u)\geq 3$
  & $Tet[t(u),n;3]$\\
  & & & $\{n,t(u)\}\not=\{5,3\}$ & \\
4 & \rule[-2ex]{0ex}{6ex} $-2$
  & $2\cos(2\pi/m)$, $m\geq 5$, $(m,2)=1$
  & $\gamma^2+4\gamma$ & $Tet[4,m;3]$\\
5 & \rule[-2ex]{0ex}{4ex} $-3$ & $(\sqrt{5}-1)/2$
  & $\sqrt{5}-1$ & $Tet[4,5;3]$\\
6 & \rule[-2ex]{0ex}{4ex} $-3$
  & $2\cos(2\pi/q)$, $q\geq 7$, $(q,4)=1$
  & $2\gamma$ & $Tet[3,4;q]$\\
7 & \rule[-2ex]{0ex}{4ex} $-3$
  & $(\sqrt{5}-3)/2$
  & $2(7+3\sqrt{5})\cosh^2u-4$, $t(u)\geq 3$
  & $Tet[3,t(u);5]$\\
8 & \rule[-2ex]{0ex}{4ex}$(\sqrt{5}-5)/2$ & $(\sqrt{5}-1)/2$
  & $(3\sqrt{5}-1)/2$ & $Tet[3,3;5]$\\
9 & \rule[-2ex]{0ex}{4ex} $-3$
   & $2\cos(\pi/m)-1$, $m\geq 4$, $(m,3)=1$
   & $\gamma^2+4\gamma$
   & $Tet[2,3,3;2,3,m]$\\
10 & \rule[-2ex]{0ex}{4ex} $-4\sin^2\frac{\pi}n$, $n\geq 5$, $(n,6)=1$
   & $\beta+3$
   & $\displaystyle\frac2{\beta}\left((\beta-3)\cos\frac{\pi}n-2\beta-3\right)$
   & $H[2;3,n;2]$ \\
11 & \rule[-1ex]{0ex}{4ex}$(\sqrt{5}-5)/2$ & $(\sqrt{5}\pm1)/2$
   & $3(\sqrt{5}+1)/2$ & $H[2;5,2;3]$\\
12 & \rule[-1ex]{0ex}{4ex} $-3$ & $(\sqrt{5}\pm1)/2$
   & $\sqrt{5}$ & $H[2;3,2;5]$\\
13 & \rule[-1ex]{0ex}{4ex}$(\sqrt{5}-5)/2$ & $(\sqrt{5}-1)/2$
   & $\sqrt{5}$ & $H[2;3,2;5]$\\
14 & \rule[-1ex]{0ex}{4ex}$(\sqrt{5}-5)/2$ & $\sqrt{5}+2$
   & $(5\sqrt{5}+9)/2$ & $H[2;3,2;5]$\\
15 & \rule[-2ex]{0ex}{5ex} $-3$
   & $2\cos(2\pi/q)$, $q\geq 8$, $(q,4)=2$
   & $2\gamma$
   & $H[2;3,3;q]$\\
\hline
\end{tabular}}
\end{sidewaystable}

\addtocounter{table}{-1}
\begin{sidewaystable}
\centering
\caption{(continued)}
{\footnotesize
\begin{tabular}{rcccc}
\hline
\rule[-2ex]{0ex}{6ex}  & $\beta=\beta(f)$ & $\gamma=\gamma(f,g)$ & $\beta'=\beta(g)$ & $\Gamma=\langle f,g\rangle$ \\
\hline
16 & \rule[-2ex]{0ex}{6ex} $-4\sin^2\frac{\pi}n$, $n\geq 5$, $(n,6)=1$
   & $2(\beta+3)$
   & $-\displaystyle\frac{6}{\beta}\left(2\cos\frac{\pi}n+\beta+2\right)$
   & $R[n,2;2]$\\
17 & \rule[-2ex]{0ex}{6ex} $-4\sin^2\frac{\pi}n$,  $(n,2)=2$,
   & $4\cosh^2u+\beta$, $(t(u),2)=2$,
   & $\displaystyle \frac{4}\gamma\cosh^2v-\frac {4\gamma}\beta$,
     $t(v)\geq 3$, $(t(v),2)=1$
   & $PH[n,t(u)/2,t(v)]$\\
   & $4\leq n<\infty$ & $1/n+1/t(u)<1/2$ & & \\
18 & \rule[-2ex]{0ex}{6ex} $-4\sin^2\frac{\pi}n$,  $(n,2)=2$,
   & $4\cosh^2u+\beta$, $(t(u),2)=2$,
   & $\displaystyle \frac{4}\gamma\cosh^2v-\frac {4\gamma}\beta$,
     $t(v)\geq 4$, $(t(v),2)=2$
   & ${\cal S}_2[n,t(u)/2,t(v)/2]$\\
   & $4\leq n<\infty$ & $1/n+1/t(u)<1/2$ & & \\
19 & \rule[-2ex]{0ex}{6ex} $-4\sin^2\frac{\pi}n$, $(n,2)=2$,
   & $4\cosh^2u+\beta$, $(t(u),2)=1$,
   & $\displaystyle\frac{4(\gamma-\beta)}\gamma\cosh^2v-\frac{4\gamma}\beta$,
   & $P[n,t(u),t(v)]$\\
   & $4\leq n<\infty$ & $1/n+1/t(u)<1/2$ & $t(v)\geq 3$, $(t(v),2)=1$ & \\
20 & \rule[-2ex]{0ex}{6ex} $-4\sin^2\frac{\pi}n$, $(n,2)=2$,
   & $4\cosh^2u+\beta$, $(t(u),2)=1$,
   & $\displaystyle\frac{4(\gamma-\beta)}\gamma\cosh^2v-\frac{4\gamma}\beta$,
   & $GTet_1[n,t(u),t(v)/2]$\\
   & $4\leq n<\infty$ & $1/n+1/t(u)<1/2$ & $t(v)\geq 4$, $(t(v),2)=2$ & \\
21 & \rule[-2ex]{0ex}{6ex} $-4\sin^2\frac{\pi}n$,  $(n,2)=1$,
   & $4\cosh^2u+\beta$, $(t(u),2)=2$,
   & $\displaystyle\frac2\gamma\left(\cosh v-\cos\frac{\pi}{n}\right)
      -\frac2{\gamma\beta}\left((\gamma-\beta)^2\cos\frac{\pi}{n}+\gamma(\gamma+\beta)\right)$
   & ${\cal S}_3[n,t(u)/2,t(v)]$\\
   & $n\geq 3$ & $1/n+1/t(u)<1/2$ & & \\
22 & \rule[-2ex]{0ex}{6ex} $-3$
   & $2\cos(2\pi/n)-1$, $(n,2)=1$,
   & $\displaystyle{\frac2{\gamma}(\gamma^2+2\gamma+2)}$
   & $GTet_1[n,3,2]$\\
   & & $n\geq 7$ & & \\
23 & \rule[-2ex]{0ex}{6ex} $-4\sin^2\frac{\pi}n$,  $(n,2)=1$,
   & $4\cosh^2u+\beta$, $(t(u),2)=1$,
   & $\displaystyle\frac{2(\gamma-\beta)}{\gamma}\cosh v
      -\frac2{\gamma\beta}\left((\gamma-\beta)^2\cos\frac{\pi}{n}+\gamma(\gamma+\beta)\right)$
   & $GTet_2[n,t(u),t(v)]$\\
   & $n\geq 3$ & $1/n+1/t(u)<1/2$ & & \\
24 & \rule[-2.5ex]{0ex}{7ex} $-4\sin^2\frac{\pi}n$, $(n,2)=1$,
   & $(\beta+4)(\beta+1)$
   & $\displaystyle\frac{2(\beta+2)^2}{\beta+1}\left(\cosh v-\cos\frac{\pi}{n}\right)
      -\frac 2\beta\left(\beta^2+6\beta+4\right)$
   & $GTet_2[n,3,t(v)]$\\
   & $n\geq 7$ & & & \\
\hline
\end{tabular}}
\end{sidewaystable}

\section{Geometry of the class $\cal D$: Kleinian orbifolds}

In this section we describe all Kleinian orbifolds $Q(\Gamma)$
for which $\Gamma\in{\cal D}$.

\begin{figure}[htbp]
\centering
\begin{tabular}{ccc}
\hline
\includegraphics[width=3 cm]{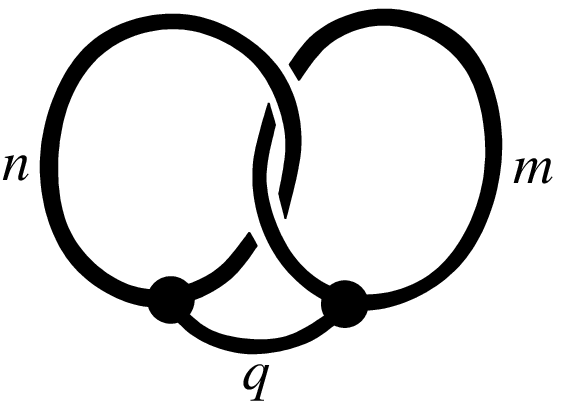}\qquad &
\qquad\includegraphics[width=2.3 cm]{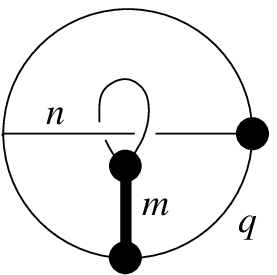} &
\qquad \includegraphics[width=2.5 cm]{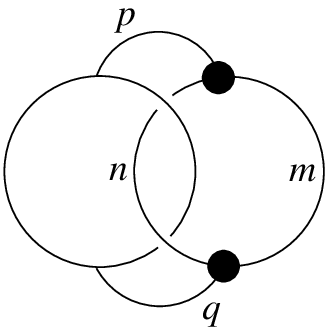}\\
(a) $GT[n,m;q]$\qquad &
\qquad(b) $PH[n,m,q]$ \qquad &
\qquad (c) $H[p;n,m;q]$
\end{tabular}

\begin{tabular}{cc}
\includegraphics[width=2.3 cm]{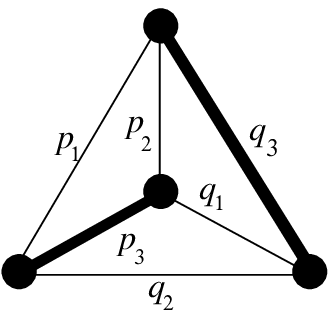}\qquad &
\qquad\includegraphics[width=2.3 cm]{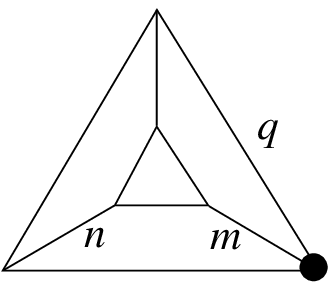}\\
(d) $Tet[p_1,p_2,p_3;q_1,q_2,q_3]$ \qquad &
\qquad (e) $P[n,m,q]$\\
\multicolumn{2}{c}{\rule[-0ex]{0ex}{5ex}Part I: Orbifolds
$Q$ embedded in ${\mathbb S}^3$ and $\pi_1^{orb}(Q)$}
\end{tabular}

\begin{tabular}{cc}
\hline
\includegraphics[width=4.5 cm]{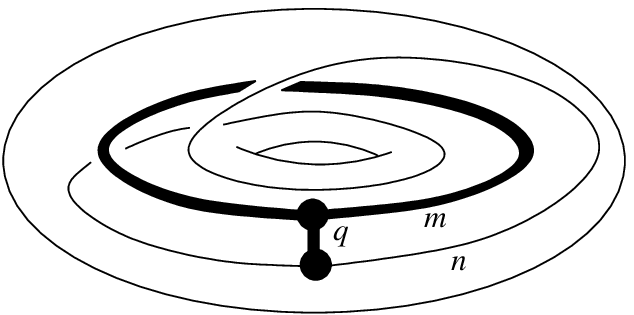}&
\includegraphics[width=4.5 cm]{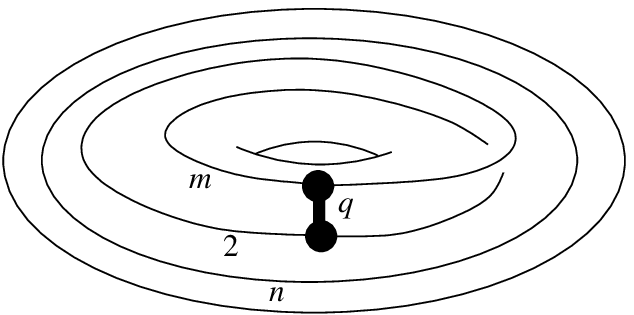}\\
(f) Orbifolds embedded in ${\cal S}(2)$; &
(g) Orbifolds embedded in ${\cal S}(2)$;\\
$\pi_1^{orb}(Q)\cong {\cal S}_2[n,m,q]$ &
$\pi_1^{orb}(Q)\cong GTet_2[n,m,q]$\\
\includegraphics[width=4.5 cm]{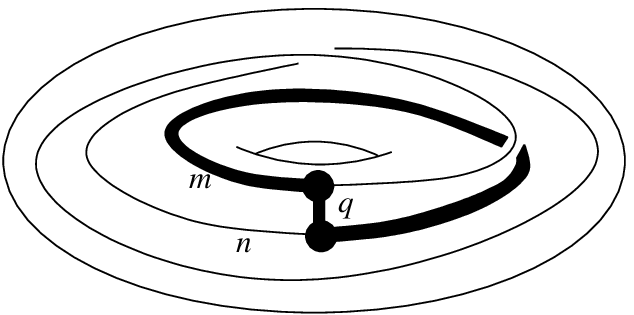}&
\includegraphics[width=4.5 cm]{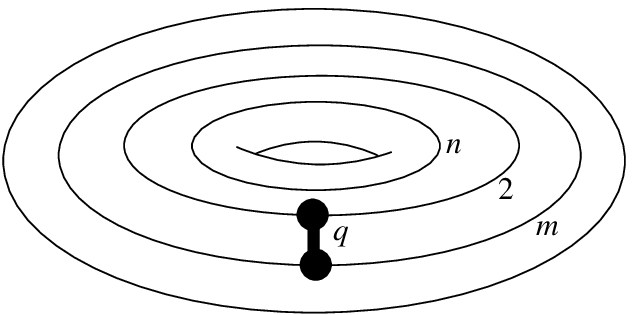}
\\
(h) Orbifolds embedded in ${\cal S}(3)$; &
(i) Orbifolds embedded in ${\mathbb S}^2\times {\mathbb S}^1$;\\
$\pi_1^{orb}(Q)\cong {\cal S}_3[n,m,q]$ &
$\pi_1^{orb}(Q)\cong GTet_1[n,m,q]$\\
\multicolumn{2}{l}{\rule[-0ex]{0ex}{5ex}Part II: Orbifolds 
$Q$ embedded
in Seifert fibre spaces; only the torus that}\\
\multicolumn{2}{l}{contains the singular set or boundary components
is shown}\\
\hline
\multicolumn{2}{c}{\includegraphics[width=2.5 cm]{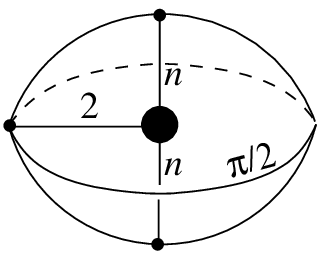}}\\
\multicolumn{2}{c}{(j) $\pi_1^{orb}(Q)\cong R[n,2;2]$}\\
\multicolumn{2}{c}{\rule[-0ex]{0ex}{5ex}Part III: Orbifolds 
$Q$ embedded in
$\mathbb{R}P^3$}\\
\hline
\end{tabular}

\caption{All orbifolds $Q(\Gamma)$ with $\Gamma\in{\cal
D}$}\label{orbifolds}
\end{figure}

\begin{theorem}\label{orb_th}
All Kleinian orbifolds $Q(\Gamma)$ with $\Gamma\in{\cal D}$ are listed
in Figure \ref{orbifolds}. The admissible
values of $n$, $m$, $p$, $q$, $p_i$, and $q_i$ 
are specified in
Theorem~\ref{allgroups}.
\end{theorem}

\begin{remark}
{\rm 
Notice that
if the orbifolds in Figure~\ref{orbifolds} 
are labelled  not like described in Theorem~\ref{allgroups},
this does not mean that the orbifolds are not hyperbolic. Sometimes, 
an orbifold
remains hyperbolic, though its group does not belong to the class
$\cal D$ anymore.

For example, $Q=Q(R[n,m;q])$ is hyperbolic for all $n,m,q\in{\mathbb Z}$
for which the determinant of the Gram matrix $\Delta$ is negative, where
$$
\Delta=\left[
\begin{array}{cccc}
1 & -\cos(\pi/q) & 0 & -\cos(\pi/2m)\\
-\cos(\pi/q) & 1 & -\cos(\pi/2m) & 0\\
0 & -\cos(\pi/2m) & 1 & -\cos(\pi/n)\\
-\cos(\pi/2m) & 0 & -\cos(\pi/n) & 1
\end{array}
\right].
$$
}
\end{remark}

\subsection{Fat edges and fat vertices}

In figures, we schematically draw singular sets and boundary components
of orbifolds using fat vertices and fat edges. In fact, each picture
gives rise to an infinite series of orbifolds which might be
compact or non-compact of finite or infinite volume.

We say that a finite 3-regular graph $\Sigma(Q)$ 
with fat vertices and fat edges
{\it represents the singular set and/or
boundary components of~$Q=Q(\Gamma)$} if
\begin{itemize}
\item[1.] edges of $\Sigma(Q)$ are labelled by positive integers 
greater than 1
or symbols $\infty$ and $\overline\infty$;
\item[2.] the endpoints of a fat edge are fat vertices;
\item[3.] if $p$, $q$, and $r$ are labels on the edges incident
to a non-fat vertex, then $1/p+1/q+1/r>1$.
\end{itemize}

To reproduce the orbifold $Q$ from a graph $\Sigma(Q)$
we first work out all fat vertices and then all fat edges
as follows.

Let $v\in\Sigma(Q)$ be a fat vertex and $p$, $q$, and
$r$ be the indices at the edges incident to~$v$.

Suppose that all $p,q,r<\infty$.
If $1/p+1/q+1/r>1$ then the vertex $v$ is a singular point of $Q$
and the local group of $v$ is one of the finite groups $D_{2n}$,
$S_4$, $A_4$, $A_5$.
If $1/p+1/q+1/r=1$ then $v$
represents a puncture.
A cusp neighbourhood of $v$ is a quotient of a horoball in
${\mathbb H}^3$ by a Euclidean triangle group $(2,3,6)$,
$(2,4,4)$, or $(3,3,3)$.
In case $1/p+1/q+1/r<1$ the vertex $v$
must be removed together with its open neighbourhood, which
means that $Q$ has a boundary component.

If one of the indices, say $p$, equals $\infty$ and
$1/p+1/q+1/r=1$, then $q=r=2$ and $v$ is a puncture.

\begin{figure}[htbp]
\centering
\begin{tabular}{ccc}
\includegraphics[width=2 cm]{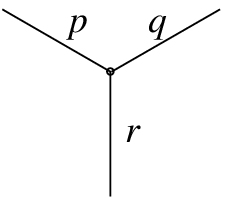} \qquad&
\qquad\includegraphics[width=2 cm]{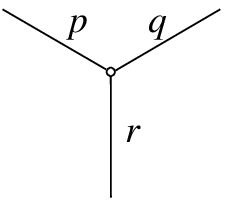} \qquad&
\qquad\includegraphics[width=2 cm]{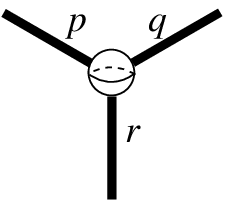} \\
$\frac1p+\frac1q+\frac1r>1$ \qquad&
\qquad$\frac1p+\frac1q+\frac1r=1$ \qquad&
\qquad$\frac1p+\frac1q+\frac1r<1$ or \\
$p,q,r<\infty$ \qquad & \qquad$p,q,r\leq\infty$ \qquad & 
\qquad$p=q=2$, $r=\overline\infty$
\end{tabular}
\caption{Fat vertices}\label{fatvert}
\end{figure}

For all the other $p$, $q$, $r$, the vertex $v$ is removed
together with its open neighbourhood.
Thus, we have three possible configurations for a fat vertex~$v$,
which are shown in Figure~\ref{fatvert}.

\medskip

Now we proceed with the fat edges.
If an edge $e$ is labelled by an integer $p<\infty$,
then $e$ is a part of the singular set of the
orbifold $Q$ and consists of cone points of index~$p$.

\begin{figure}[htbp]
\begin{center}
\begin{tabular}{cc}
\includegraphics[width=2 cm]{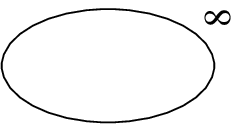}\quad&
\quad\includegraphics[width=6 cm]{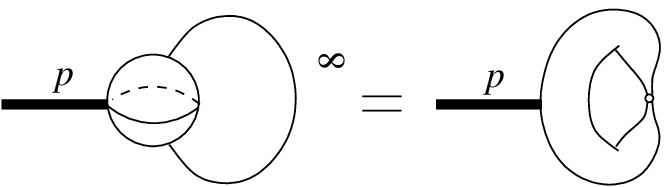}\\
non-rigid rank 2 cusp;  \quad & 
\quad rank 1 cusp; $F$ is an annulus\\
$F$ is a torus \quad & \\
(a) \quad & \quad (b) \rule[-6ex]{0ex}{8ex}
\end{tabular}
\begin{tabular}{ccc}
\includegraphics[width=3.5 cm]{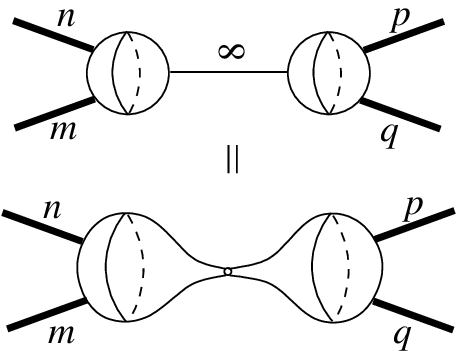}\quad&
\quad\includegraphics[width=3 cm]{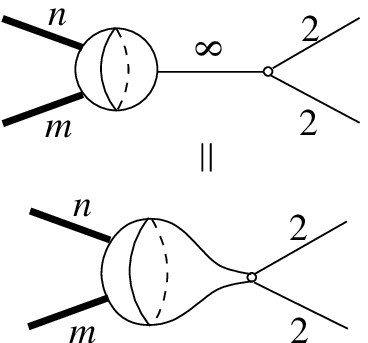}\quad&
\quad\includegraphics[width=2.5 cm]{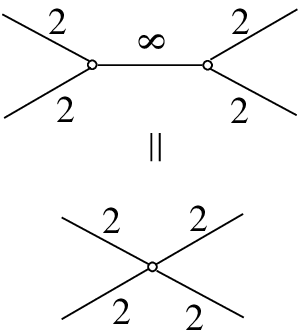}\\
rank 1 cusp; \quad & \quad rank 1 cusp; $F$ is a 
\quad & \quad non-rigid rank 2 cusp;\\
$F$ is an annulus \quad & \quad disc with two cone \quad & 
\quad $F$ is a pillow\\
& \quad points of order 2 \quad & \\
(c) \quad & \quad (d) \quad & \quad (e)\\
\end{tabular}
\end{center}
\caption{Cusps}\label{cusps}
\end{figure}

Edges marked by $\infty$ represent cusps of~$Q$.
A cusp neighbourhood is the quotient of a horoball by an
elementary parabolic group.
Topologically it is $F\times [0,\infty)$, where $F$ is a Euclidean
orbifold called the cross-section of the cusp
(see, e.g., \cite{BMP03} for geometric structures on orbifolds).
All cusps are drawn in Figure~\ref{cusps}. Note that for
the orbifolds $Q(\Gamma)$ with $\Gamma\in\cal D$, the cases
shown in Figures~\ref{cusps}(a) and~(e) never occur.

If $e$ is labelled by $\overline\infty$, then it must be deleted
together with its open regular neighbourhood.

\subsection{Geometry of the orbifolds}

\subsubsection*{Orbifolds embedded in ${\mathbb S}^3$ or 
$\mathbb{R}P^3$}

The orbifolds $Q$ shown in Figures~\ref{orbifolds}(a)--\ref{orbifolds}(e)
are embedded in ${\mathbb S}^3={\mathbb R}^3\cup\{\infty\}$ 
so that $\infty$ is a
non-singular interior point of~$Q$.

In Figure~\ref{orbifolds}(j), $\mathbb{R}P^3$ is shown as a lens with
antipodal points on the boundary identified. The angle at 
the edge of the lens is $\pi/2$ and, therefore, the edge is mapped
onto a singular loop with index~2.

\subsubsection*{Orbifolds embedded in ${\cal S}(n)$ or 
${\mathbb S}^2\times {\mathbb S}^1$}

Let $T(n)$ be a Seifert fibred solid torus obtained from a trivial
fibred solid torus $D^2\times {\mathbb S}^1$ by cutting it along
$D^2\times \{x\}$ for some $x\in {\mathbb S}^1$, rotating one of the
discs through $2\pi/n$ and gluing back together.

Denote by
${\cal S}(n)$ a space obtained by gluing two copies of $T(n)$
along their boundaries fibre to fibre.
Clearly, ${\cal S}(n)$ is homeomorphic to 
${\mathbb S}^2\times {\mathbb S}^1$
and is $n$-fold covered by trivially fibred 
${\mathbb S}^2\times {\mathbb S}^1$.
There are two critical fibres whose length is $n$ times shorter
than the length of a regular fibre.

The orbifolds shown in Figures~\ref{orbifolds}(f)--\ref{orbifolds}(i)
are embedded in Seifert fibre spaces ${\cal S}(n)=T(n)\cup T(n)$
or ${\mathbb S}^2\times {\mathbb S}^1$ trivially fibred.
We draw only the solid torus that contains
singular points (or boundary components). The other fibred torus is
meant to be attached and is not shown.

For example, if $\Gamma={\cal S}_2[n,m,q]$, then $Q(\Gamma)$
is embedded into ${\cal S}(2)$
and $\Sigma(Q)$ is placed in ${\cal S}(2)$ so
that the edge of order $m$ (if $m<\infty$)
lies on a critical fibre
of ${\cal S}(2)$ and the edge of order $n$ lies on a regular one.
Note that the generator $L$ of the group ${\cal S}_2[n,m,q]$
is $\pi$-loxodromic~\cite{KK05}.

\subsection{Orbifolds of finite volume}

In this section, we list all finite volume orbifolds
$Q(\Gamma)$ for $\Gamma\in\D$.

\subsubsection*{Compact orbifolds}

\begin{itemize}
\item[(b)] $\pi_1^{orb}(Q)\cong PH[n,m,q]$: $n=4$, $3\leq m\leq 5$, $q=3$
\item[(c)] $\pi_1^{orb}(Q)\cong H[p;n,m;q]$:
$[p;n,m;q]$ is $[2;2,3;5]$, or $[2;2,5;3]$, or $[2;3,5;2]$
\item[(d)] $\pi_1^{orb}(Q)\cong Tet[2,3,3;2,3,q_3]$: $q_3=4,5$

$\pi_1^{orb}(Q)\cong Tet[n,m;q]$: $n=5$, $m=4,5$, $q=3$ or
$n=m=3$, $q=5$
\item[(e)] $\pi_1^{orb}(Q)\cong P[n,m,q]$:

$8\leq n<\infty$, $n$ is even, $m=3$, $q=3,5$; 

$4\leq n<\infty$, $n$ is even, $m=5$, $q=3$
\item[(g)] $\pi_1^{orb}(Q)\cong GTet_2[n,m,q]$:

$n\geq 7$ is odd, $m=3$, $3\leq q\leq 5$; 

$n\geq 5$ is odd, $m=5$, $q=3$;

$n,m\geq 3$ are odd, $1/n+1/m<1/2$, $q=2$
\item[(h)] $\pi_1^{orb}(Q)\cong {\cal S}_3[n,m,q]$:

$n\geq 5$ is odd, $m=q=2$;

$n=5$, $m=2$, $q=3$;

$n=5$, $m=3$, $q=2$;

$n=3$, $m=4,5$, $q=2$
\item[(i)] $\pi_1^{orb}(Q)\cong GTet_1[n,m,q]$:
$7\leq n<\infty$, $m=3$, $q=2$ 
\end{itemize}

There are no compact orbifolds in the class $\cal D$ with
$\pi_1^{orb}(Q)\cong GT[n,m;q]$, ${\cal S}_2[n,m,q]$, or
$R[n,2;2]$.

\subsubsection*{Non-compact orbifolds of finite volume}

\begin{itemize}
\item[(a)] $\pi_1^{orb}(Q)\cong GT[n,m;q]$: $[n,m;q]$ is
$[3,3;3]$, or $[4,4;2]$, or $[4,3;2]$
\item[(b)] $\pi_1^{orb}(Q)\cong PH[n,m,q]$:
$n=4$, $m=6$, $q=3$ or $n=6$, $2\leq m\leq 6$, $q=3$
\item[(d)] $\pi_1^{orb}(Q)\cong Tet[n,m;q]$:
$n=6$, $3\leq m\leq 6$, $q=3$
\item[(f)] $\pi_1^{orb}(Q)\cong {\cal S}_2[n,m,q]$:
$n=4$, $m=3,4$, $q=2$
\item[(g)] $\pi_1^{orb}(Q)\cong GTet_2[n,m,q]$:
$n\geq 7$ is odd, $m=3$, $q=6$
\item[(h)] $\pi_1^{orb}(Q)\cong {\cal S}_3[n,m,q]$:
$n=3$, $m=6$, $q=2$
\item[(i)] $\pi_1^{orb}(Q)\cong GTet_1[n,m,q]$:
$8\leq n<\infty$, $(n,2)=2$, $m=3$, $q=3$
\end{itemize}

There are no non-compact orbifolds of finite volume in the class $\cal D$
with $\pi_1^{orb}(Q)\cong H[p;n,m;q]$, $P[n,m,q]$, or $R[n,2;2]$.

\appendix
\section{Both generators are elliptic}

Suppose an $\RP$ group $\Gamma$ is generated by two elliptic 
elements $f$ and $g$.
By \cite[Theorem~2]{KK02}, the axes of $f$ and $g$ either lie in one
plane or are mutually orthogonal skew lines. Note that such a group is
elementary or has an invariant plane if and only if the axes 
of $f$ and $g$ lie in
one plane or one of the generators has order~2.
So we do not consider such groups.

Theorem~\ref{two_ell} gives necessary and sufficient
discreteness conditions for
non-elementary $\RP$ groups without invariant plane generated by two
primitive elliptics. The proof of Theorem~\ref{two_ell} is constructive
and we will use it later to find parameters for $\Gamma\in\D$
generated by two primitive elliptics and describe corresponding orbifolds.

\begin{theorem}[\cite{Kli89}]\label{two_ell}
Let $f$ and $g$ be primitive elliptic elements of orders $n\geq 3$ and
$m\geq 3$,
respectively, and let their axes be mutually orthogonal skew lines.
Then:
\begin{itemize}
\item[$(1)$] there exists a unique element $h\in \PSL$ such that
$h^2=fgf^{-1}g^{-1}$
and $(hg)^2=1$, and
\item[$(2)$] $\Gamma=\langle f,g\rangle$ is discrete if and only if one of
the following holds:
\begin{itemize}
\item[$(2.i)$] $h$ is hyperbolic, parabolic, or a primitive elliptic
element of order $p$, where $\cos(\pi/p)>\sin(\pi/n)\sin(\pi/m)$;
\item[$(2.ii)$] $m=n\geq 7$, $(n,2)=1$, and $h$ is a square of a primitive
elliptic element of order $n$.
\end{itemize}
\end{itemize}
\end{theorem}

\noindent{\it Proof. }
{\bf 1.} {\it Construction of $\Gamma^*$.}
We start with construction of a reflection group $\Gamma^*$ 
containing $\Gamma$
as a subgroup of finite index. Such a group is discrete if and
only if so is $\Gamma$.

Let $f$ and $g$ be primitive elliptic elements of orders $n\geq 3$
and $m\geq 3$, respectively, and let the axes of $f$ and $g$ be mutually
orthogonal skew lines.
We will denote elements and their axes by the same
letters when it does not lead to any confusion.

\begin{figure}[htbp]
\centering
\begin{tabular}{cc}
\includegraphics[width=4 cm]{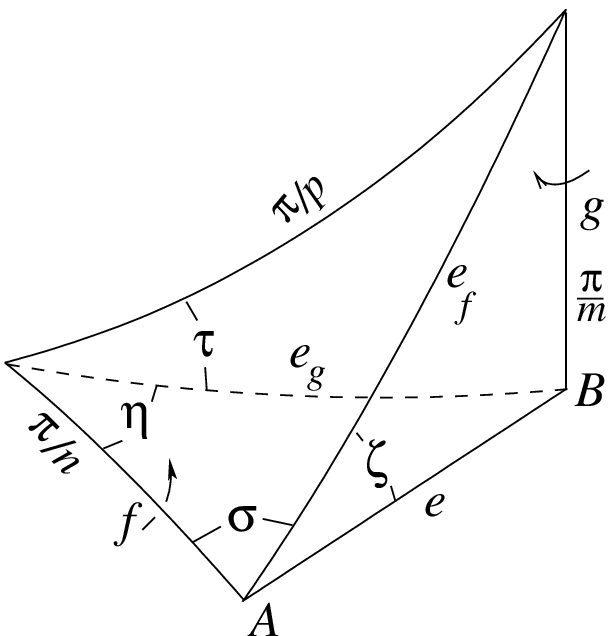} \quad & \quad
\includegraphics[width=4.5 cm]{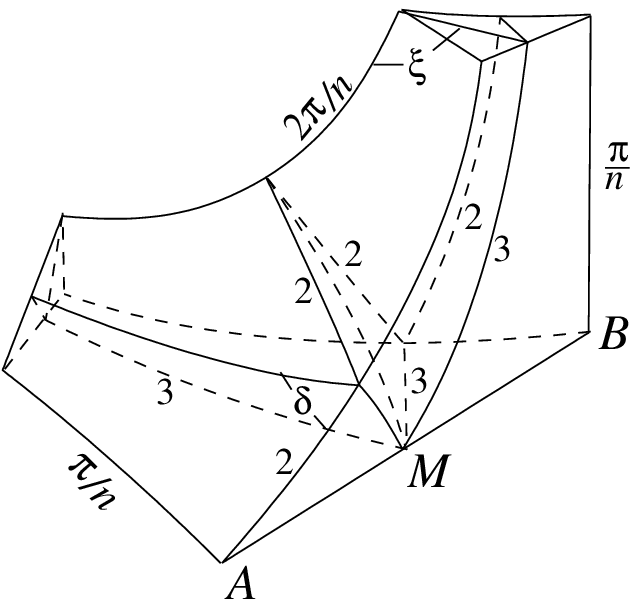}\\
(a) \quad & \quad (b)\\
\end{tabular}
\caption{}\label{fund_poly}
\end{figure}

Let $e$ be a half-turn whose axis is orthogonal to both
$f$ and $g$.
Let $e_f$ and $e_g$ be half-turns such that $f=e_fe$ and~$g=e_ge$.
Then $e_f$ and $e$ intersect at an angle of~$\pi/n$
and $e_g$ and $e$ intersect at an angle of~$\pi/m$.

Denote by $\eta$ the plane through $e_g$, $e$, and $f$ and  by $\zeta$
the plane through $e_f$, $e$, and $g$. 
Denote by $\sigma$ the plane through
$e_f$ and $f$ and by $\tau$ the plane through $e_g$ and~$g$
(see Figure~\ref{fund_poly}(a)).

For the group $\Gamma=\langle f,g \rangle$ we define two finite
index extensions of it as follows:
$\widetilde{\Gamma}=\langle f,g,e\rangle$ and
$\Gamma^*=\langle f,g,e,R_\eta\rangle$
(we denote the reflection in a plane $\kappa$ by $R_\kappa$).

It is easy to see that $\widetilde\Gamma=\Gamma\cup\Gamma e$.
If $e\in\Gamma$ then $\widetilde\Gamma=\Gamma$ and if
$e\notin\Gamma$ then $\Gamma$ is a subgroup of index~2 in
$\widetilde\Gamma$. (As we will see in Corollary~\ref{cor2}, both possibilities are
realised.) Since, moreover, $\widetilde\Gamma$ is the orientation
preserving subgroup of index~2 in $\Gamma^*$,
the groups $\Gamma$, $\widetilde\Gamma$,
and $\Gamma^*$ are either all discrete or all non-discrete.
It is clear that
$\Gamma^*=\langle R_\eta,R_\zeta,R_\sigma,R_\tau\rangle$.

\medskip\noindent
{\bf 2.}
{\it Existence and uniqueness of $h$.}
First, we show that if $h=R_\sigma R_\tau$ then $h^2=[f,g]$
and $(hg)^2=1$. Indeed,
\begin{eqnarray*}
fgf^{-1}g^{-1}&=&(e_fe)(e_ge)(ee_f)(ee_g)=(e_fee_g)^2=(fe_g)^2\\
&=&(R_\sigma R_\eta R_\eta R_\tau)^2=(R_\sigma R_\tau)^2=h^2.
\end{eqnarray*}
Moreover, $hg=(R_\sigma R_\tau)(R_\tau R_\zeta)=R_\sigma R_\zeta=e_f$.
So, $(hg)^2=e_f^2=1$. Hence, the element $h$ required in 
assertion~(1) of Theorem~\ref{two_ell} exists. 

Now show that $h$ is unique. If $[f,g]$ is parabolic, then the
square root of $[f,g]$ is unique in $\PSL$ and we are done. 
Let $[f,g]$ be elliptic or hyperbolic
(it cannot be strictly loxodromic, because 
$[f,g]=(R_\sigma R_\tau)^2=R_\sigma(R_\tau R_\sigma R_\tau)$
is the product of two reflections). Denote by $\overline h$
the second square root of $[f,g]$.
If $[f,g]$ is hyperbolic,
then $\overline h$ is $\pi$-loxodromic with the same axis and translation
length as $h$. If $[f,g]$ is elliptic, then $\overline h$ is elliptic
with the same axis as $h$, but the rotation angle of
$\overline h$ is $(\pi-2\pi/p)$.
It is clear that in both cases
$\overline h^2=[f,g]$ and $(\overline hg)^2\not=1$. Uniqueness is
proved.

\medskip\noindent
{\bf 3.} {\it Polyhedron $\cal T$.}
We first consider the polyhedron $\cal T$ bounded by the planes
$\eta$, $\zeta$, $\sigma$, and $\tau$ 
(we know that reflections in these planes generate $\Gamma^*$)
and determine the dihedral angles of $\cal T$.
The line $f$ is the line of intersection of the planes
$\eta$ and $\sigma$ and the dihedral angle made by $\eta$ and
$\sigma$ equals $\pi/n$. Similarly, $g$ is the line of intersection of
$\zeta$ and $\tau$ and the dihedral angle made by them equals $\pi/m$.
The plane $\eta$ is orthogonal to both $\zeta$ and $\tau$, and
the plane $\sigma$ is orthogonal to $\zeta$.
As for $\sigma$ and $\tau$, they may either intersect, or be parallel
or disjoint.
Denote the dihedral angle of $\cal T$ made by $\sigma$ and $\tau$ by
$\pi/p$
($p$ is not necessarily an integer).
We keep the notation $\pi/p$
with $p=\infty$ and $p=\overline\infty$
for the cases of parallel and disjoint
$\sigma$ and $\tau$, respectively.
From \cite{Vin85}, ${\cal T}$ exists in ${\mathbb H}^3$ if
and only if
\begin{equation}\label{existence_cond}
\cos(\pi/p)>\sin(\pi/n)\sin(\pi/m)
\end{equation}

Note that (\ref{existence_cond}) implies $\pi/p<\pi/2$.
Therefore, there is an obvious connection between $h=R_\sigma R_\tau$ and
the dihedral angle of $\cal T$ made by $\sigma$ and $\tau$.
In particular, $h$ is parabolic (hyperbolic) if and only if
$\sigma$ and $\tau$ are parallel (disjoint, respectively).

\medskip\noindent
{\bf 4.} {\it Discreteness of $\Gamma$ implies (2.i) or (2.ii).}

Assume that $\Gamma$ is discrete, then $\Gamma^*$ is also discrete 
and, therefore, $h=R_\sigma R_\tau$ either satisfies $(2.i)$ or
is a non-primitive elliptic element of finite order,
i.e., $p=q/k$, where $q$ and $k$ are coprime integers, $k>1$.
Let us show that then $(2.ii)$ holds.

Consider a surface $S$ that is orthogonal to the planes $\eta$, $\sigma$,
and $\tau$. If $\eta$, $\sigma$, and $\tau$ intersect at a point
$x\in{\mathbb H}^3$, then $S$ is a sphere with center $x$; if $\eta$,
$\sigma$,
and $\tau$ meet at a point $x\in\partial {\mathbb H}^3$, then $S$ is a
horosphere; if $\eta$, $\sigma$, and $\tau$ do not have a common point
in $\overline{{\mathbb H}^3}$,
then $S$ is a hyperbolic plane. Consider the intersection $\Delta$ of
$S$ with the polyhedron $\cal T$, that is a link made by $\eta$,
$\sigma$,
and $\tau$ (cf. \cite{KK02}).

The link $\Delta$ is then a spherical, Euclidean, or hyperbolic triangle
(depending on the type of $S$) with angles of $\pi/2$, $\pi/n$, and
$\pi/p=k\pi/q$,
where $k$ and $q$ are coprime and

\begin{itemize}
\item[(a)] Suppose that $\Delta$ is spherical. From the list of all
spherical triangles
with two primitive and one non-primitive angles~\cite{Fel99},
there are two possibilities for $\Delta$: $(\pi/2,\pi/5,2\pi/5)$ and
$(\pi/2,\pi/3,2\pi/5)$. Since in both cases $k\pi/q$=$2\pi/5$, the link of
the vertex made by $\zeta$, $\sigma$, and $\tau$ is also spherical,
that is $n=3$ or~$5$, and $m=3$ or~$5$. However, all triples $(n,m,p)$,
where $n,m\in\{3,5\}$ and $p=5/2$, do not satisfy~(\ref{existence_cond}).
Contradiction. So, $\Delta$ cannot be spherical.

\item[(b)] $\Delta$ also
cannot be Euclidean, because there are no Euclidean
triangles with angles $\pi/2$, $\pi/n$, $k\pi/q$ with $(k,q)=1$
that lead to a discrete group.

\item[(c)] Suppose that $\Delta$ is a hyperbolic triangle with angles
$\pi/2$, $\pi/n$, $k\pi/q$. From the list of all hyperbolic triangles
with one non-primitive and two primitive angles~\cite{Mat82}, we have
that $k=2$ and $q=n\geq 7$ is odd.

Since $q\geq 7$, the link of the vertex made by $\zeta$, $\sigma$, and
$\tau$ is also a hyperbolic triangle and again from~\cite{Mat82}
we conclude
that $m=q=n$.

\end{itemize}

From (a), (b), and (c) above, it follows that the only possibility for
dihedral angles of $\cal T$ is that with $n=m\geq 7$, $(n,2)=1$, and the
dihedral angle between $\sigma$ and $\tau$ equals $2\pi/n$, that is
$h=R_\sigma R_\tau$ is a square of a primitive elliptic element
of order $n$. Hence, $(2.ii)$ holds.

\medskip\noindent
{\bf 5.} {\it $(2.i)$ implies discreteness of $\Gamma$.}
Assume that $(2.i)$ holds. Then the inequality (\ref{existence_cond})
implies that the polyhedron $\cal T$ exists in $\mathbb H^3$. Moreover, since
$h$ is a hyperbolic, parabolic, or primitive elliptic element of
order $p$, the planes $\sigma$ and $\tau$ are either disjoint,
or parallel, or the dihedral angle between them is $\pi/p$, 
$p\in\mathbb Z$. Therefore, $\cal T$
and reflections $R_\eta$, $R_\zeta$, $R_\sigma$, and
$R_\tau$
satisfy the hypotheses of the Poincar\'e theorem~\cite{EP94},
$\Gamma^*$ is discrete,
and $\cal T$ is its fundamental polyhedron.

\medskip\noindent
{\bf 6.} {\it $(2.ii)$ implies discreteness of $\Gamma$.}
Since $h$ is a square of a primitive elliptic element
of order $n$ and $m=n$,
the dihedral angles of $\cal T$ are as described in part 3 of the proof
with $\pi/m=\pi/n$ and $\pi/p=2\pi/n$.

Note that for each $n$, ${\cal T}={\cal T}(n)$ is uniquelly determined 
by its dihedral angles. 
(To see this, we apply Andreev's thorem to the compact polyhedron
that has two more faces: the hyperbolic plane $S$ orthogonal
to $\eta$, $\sigma$, and $\tau$, and the plane $\widetilde S$
orthogonal to $\zeta$, $\sigma$, $\tau$.) 
In particular, the length of $AB$ in Figure~\ref{fund_poly}(b)
is determined by the dihedral angles of $\cal T$ uniquely, and
so it depends only on~$n$.

Since the dihedral angle between $\sigma$ and $\tau$ is $2\pi/n$, there is
one more reflection $R_\xi$ in $\Gamma^*$, 
$R_\xi\in\langle R_\sigma, R_\tau\rangle$, where $\xi$ is the bisector of
the dihedral angle of $\cal T$ made by $\sigma$ and $\tau$. 
The plane $\xi$ divides $\cal T$ into
two symmetric polyhedra, each of which is further subdivided
(by planes of reflections from $\Gamma^*$) into three tetrahedra
$T(n)=T[2,2,3;2,3,n]$ with primitive dihedral angles.
To see this, look at the tesselation of $\mathbb H^3$ by the 
infinite volume
tetrahedra $T(n)$ as above and find the polyhedron ${\cal T}(n)$ which is
formed by six tetrahedra~$T(n)$.

On the other hand, reflections in the faces of ${\cal T}(n)$
generate all reflections in the faces of~$T(n)$. Since ${\cal T}(n)$
is determined by its dihedral angles uniquely, we conclude that
$\Gamma^*$ is the group of reflections in the faces
of $T(n)$ and, therefore, both $\Gamma^*$ and $\Gamma$ are discrete.
\qed

The following corollary is just a reformulation of Theorem~\ref{two_ell}
in terms of parameters.

\begin{cor}\label{cor1}
Let $f,g\in\PSL$, $\beta(f)=-4\sin^2(\pi/n)$,
$\beta(g)=-4\sin^2(\pi/m)$, where $n,m\in{\mathbb Z}$ and $n,m\geq 3$,
and let $\gamma(f,g)<-\beta(f)\beta(g)/4$.
Then $\Gamma=\langle f,g\rangle$
is discrete if and only if one of the following holds:
\begin{enumerate}
\item $\gamma(f,g)\in(-\infty;-4]$;
\item $\gamma(f,g)=-4\cos^2(\pi/p)$, $p\in{\mathbb Z}$,
$\cos(\pi/p)>\sin(\pi/n)\sin(\pi/m)$;
\item $m=n\geq 7$, $(n,2)=1$, $\gamma(f,g)=-(\beta+2)^2$.
\end{enumerate}
\end{cor}

\noindent
{\it Proof.}  Since $\beta(f)=-4\sin^2(\pi/n)$,
$\beta(g)=-4\sin^2(\pi/m)$, $n,m\in{\mathbb Z}$,
$f$ and $g$ are primitive elliptic elements. Since
$\gamma(f,g)<-\beta(f)\beta(g)/4$, the axes of $f$ and $g$ are
mutually orthogonal skew lines \cite[Theorem~4]{KK02}. So the
hypotheses of Corollary~\ref{cor1} are equivalent to the hypotheses of
Theorem~\ref{two_ell}. 
Therefore, in order to prove Corollary~\ref{cor1} it suffices
to compute parameters for each discrete group described in
Theorem~\ref{two_ell}.

We need to rewrite
conditions $(2.i)$ and $(2.ii)$ via $\gamma(f,g)$.
Since $\gamma(f,g)=\tr[f,g]-2$ and $h$ is a square root of $[f,g]$,
it is not difficult to find $\gamma(f,g)$.

The element $h$ is hyperbolic if and only if the planes
$\sigma$ and $\tau$ are disjoint.
(We denote planes (and lines) as in the proof of Theorem~\ref{two_ell}.)
Let $d$ be the hyperbolic distance between them. Since
$[f,g]=h^2=(R_\sigma R_\tau)^2$,
$$
\gamma(f,g)=\tr[f,g]-2=-2\cosh(2d)-2<-4
$$
(we must take $\tr[f,g]$ to be negative, because $\gamma(f,g)$
is strictly negative for all values $\beta(f)$, $\beta(g)$).

The element $h$ is parabolic if and only if $[f,g]$ is parabolic and
if and only if $\tr[f,g]=-2$, that is,
$\gamma(f,g)=\tr[f,g]-2=-4$.

Thus, $h$ is hyperbolic or parabolic if and only if
$\gamma(f,g)\in(-\infty,-4]$, and part~1 of Corollary~\ref{cor1}
is proved.

\medskip
Now suppose that $h$ is an elliptic element with rotation angle $\phi$,
where $\phi/2=\pi/p<\pi/2$ is the dihedral angle of $\cal T$ between
$\sigma$ and $\tau$. Then $[f,g]=h^2$ is also elliptic with rotation angle
$2\phi$. Since $\tr[f,g]$ is well-defined (does not depend on the choice
of representatives of $f$ and $g$ in ${\rm SL}(2,{\mathbb C})$) we can
determine which formula, $\tr[f,g]=+2\cos\phi$
or $\tr[f,g]=-2\cos\phi$, is correct. The easiest way to do this
is by using continuity of $\tr[f,g]$ as a function of $\phi$ and
the limit condition $\tr[f,g]\to -2$ as $\phi\to 0$.
So we must take $\tr[f,g]=-2\cos\phi$ where $\phi<\pi$ is the
doubled dihedral angle of $\cal T$.

On the other hand, if $\tr[f,g]$ is given, we can use the formula
$\tr[f,g]=-2\cos\phi$, $\phi<\pi$, to obtain the rotation angle $\phi$
of the element $h$ from Theorem~\ref{two_ell}.

Thus, $h$ is a primitive elliptic element of order $p$ if and only if
$\phi=2\pi/p$ which is equivalent to
$$
\gamma(f,g)=\tr[f,g]-2=-2\cos(2\pi/p)-2=-4\cos^2(\pi/p), \ p\in{\mathbb Z},
$$
which is part~2 of Corollary~\ref{cor1}.

At last, $h$ is a square of a primitive elliptic element of order~$n$
if and only if $\phi=4\pi/n$, which corresponds to
$\gamma(f,g)=-4\cos^2(2\pi/n)=-(\beta+2)^2$.

Corollary~\ref{cor1} is proved.
\qed

\begin{remark}\label{rem_funcmu}
{\rm
The parts 1 and 2 of Corollary~\ref{cor1} can be rewritten in terms of
the function $t(w)$ as $\gamma(f,g)=-4\cosh^2w$, where $w\in\U$
and $\cos(\pi/t(w))>\sin(\pi/n)\sin(\pi/m)$. Indeed,
\begin{eqnarray*}
\gamma=-4\cos^2(\pi/p),\ p\in{\mathbb Z}
&\quad{\Leftrightarrow}\quad&
\gamma=-4\cosh^2w,\ w\in\U,\ t(w)<\infty,\\
\gamma=-4
&\quad{\Leftrightarrow}\quad&
\gamma=-4\cosh^2w,\ w\in\U,\ t(w)=\infty,\\
\gamma<-4,
&\quad{\Leftrightarrow}\quad&
\gamma=-4\cosh^2w,\ w\in\U,\ t(w)=\overline\infty.
\end{eqnarray*}
Keep in mind that in this notation, $(t(w),2)=1$ implies that $t(w)$
is
finite and odd, while $(t(w),2)=2$ means that $t(w)$ can be even,
$\infty$, or $\overline\infty$.
}
\end{remark}

The following corollary gives group presentations for all 
groups from class $\D$ generated by two elliptics.
The corresponding orbifolds are drawn in Figures~\ref{orbifolds}(a)
and~\ref{orbifolds}(d).

\begin{cor}\label{cor2}
Let $\Gamma=\langle f,g\rangle\in{\cal D}$, $\beta(f)=-4\sin^2(\pi/n)$,
$\beta(g)=-4\sin^2(\pi/m)$, where $n,m\in{\mathbb Z}$, $n,m\geq 3$.
\begin{enumerate}
\item If $\gamma(f,g)=-4\cosh^2w$, where $w\in\U$, $(t(w),2)=2$, and
$\cos(\pi/t(w))>\sin(\pi/n)\sin(\pi/m)$,
then $\Gamma\cong GT[n,m;t(w)/2]$.
\item If $\gamma(f,g)=-4\cosh^2w$, where $w\in\U$, $(t(w),2)=1$, and
$\cos(\pi/t(w))>\sin(\pi/n)\sin(\pi/m)$,
then $\Gamma\cong Tet[n,m;t(w)]$.
\item If $m=n\geq 7$, $(n,2)=1$, and $\gamma(f,g)=-(\beta+2)^2$,
then $\Gamma\cong Tet[3,n;3]$.
\end{enumerate}
\end{cor}

\noindent
{\it Proof.} All parameters for discrete groups in the statement of
Corollary~\ref{cor2} are described in Corollary~\ref{cor1}. We will obtain
a presentation for each case by using the Poincar\'e polyhedron theorem.

We start with construction of a fundamental polyhedron and a presentation
for
the group $\widetilde\Gamma$ defined in the proof of
Theorem~\ref{two_ell}.
Since $\widetilde\Gamma$ is an orientation preserving index~2 subgroup
in $\Gamma^*$ and $\cal T$ is a fundamental polyhedron for $\Gamma^*$,
a fundamental polyhedron $\widetilde{\cal T}$ for $\widetilde\Gamma$
consists of two copies
of~$\cal T$.
We take $\widetilde{\cal T}$ to be the polyhedron bounded by
$\eta$, $\sigma$, $\tau$, and $R_\zeta(\tau)$.

By applying the Poincar\'e theorem to $\widetilde{\cal T}$ and face
pairing
transformations $e$, $e_f$, and $g$, we get
$$
\widetilde\Gamma=\langle e,e_f,g\,|\,e^2,e_f^2,g^m,(e_fe)^n,(ge)^2,
(ge_f)^p\rangle,
$$
where $p$ is an integer, $\infty$, or $\overline\infty$.
Since $f=e_fe$,
$$
\widetilde\Gamma=\langle
f,g,e\,|\,f^n,g^m,e^2,(fe)^2,(ge)^2,(gfe)^p\rangle.
$$
Note that if $p$ is odd, then using the relations
$(fe)^2=(ge)^2=1$, we see that from $(gfe)^p=1$ it follows that
$e=f^{-1}g^{-1}(fgf^{-1}g^{-1})^{(p-1)/2}$. Hence, in this case
$\widetilde\Gamma=\Gamma$ and $\Gamma\cong Tet[n,m;p]$.
Identifying faces of $\widetilde{\cal T}$, we get an orbifold
${\mathbb H}^3/\Gamma$ (see Figure~\ref{orbifolds}(d)).

If $p$ is even, $\infty$, or $\overline\infty$,
then $\Gamma$ is a subgroup of index~2 in $\widetilde\Gamma$.
To see this we apply the Poincar\'e theorem to a polyhedron
bounded by $\sigma$, $\tau$, $R_\eta(\sigma)$, and $R_\zeta(\tau)$.
Then
$$
\Gamma=\langle f,g\,|\,f^n,g^m,(fgf^{-1}g^{-1})^{p/2}\rangle.
$$
The orbifold ${\mathbb H}^3/\Gamma$ is shown in Figure~\ref{orbifolds}(a).
The parts~1 and~2  are proved.

Now suppose that we are in the case~3 of Corollary~\ref{cor1} and obtain
a presentation for $\Gamma$.

Let $\delta$ be the plane shown in Figure~\ref{fund_poly}(b) and let
$u=R_\sigma R_\delta$.
It is clear that
$$
G=\langle f,e_f,u\,|\,f^n,e_f^2,u^2,(fe_f)^2,(ue_f)^3,(fu)^3\rangle
\cong Tet[3,n;3].
$$

We claim that $G$
is isomorphic to $\Gamma=\langle f,g\rangle$.
From the link made by $\zeta$, $\sigma$, and $\tau$, we have that
$g=uf^2ue_f\in G$. Hence, $\Gamma=\langle f,g\rangle$ is a subgroup
of~$G$.

Let us show that $G$ is a subgroup of~$\Gamma$, i.e., that
$e_f$ and $u$ belong to~$\Gamma$.
From the link made by $\eta$, $\sigma$, and $\tau$, we get
$uf^2u=h^{-1}$. Since $h^n=1$ and $n$ is odd,
$$
h=(h^2)^{-\frac{n-1}{2}}=[f,g]^{-\frac{n-1}{2}}.
$$
Moreover, since $f^n=1$,
\begin{equation}\label{ufu}
ufu=(uf^{-2}u)^{\frac{n-1}2}=h^{\frac{n-1}2}=[f,g]^{-\frac{(n-1)^2}4}.
\end{equation}
Further, combining~(\ref{ufu}) with the relation $(fu)^3=1$, we have
\begin{eqnarray*}
u&=&f(ufu)f=f[f,g]^{-\frac{(n-1)^2}4}f\in\Gamma {\rm \quad and}\\
e_f&=&uf^{-2}ug=[f,g]^{\frac{(n-1)^2}2}g\in\Gamma.
\end{eqnarray*}
Thus, $G$ is a subgroup of $\Gamma$ and so $\Gamma\cong G$.
\qed


\begin{thebibliography}{99}
%
\bibitem{BP01}
M.~Boileau and J.~Porti,
{\it Geometrization of $3$-orbifolds of cyclic type. With an appendix:
Limit of hyperbolicity for spherical $3$-orbifolds by Michael
Heusener and Joan Porti},
Ast\'erisque. 272. Paris: Sosi\'et\'e Math\'ematique de France.
2001.
%
\bibitem{BMP03}
M.~Boileau, S.~Maillot, J.~Porti,
{\it Three-dimensional orbifolds and their geometric structures},
Panoramas et Synth\`eses. 15. Sosi\'et\'e Math\'ematique de France.
2003.
%
\bibitem{EP94}
D.~B.~A.~Epstein and C.~Petronio, {\it An exposition of Poincar\'e's
polyhedron theorem},
L'Enseignement Math\'ematique {\bf 40} (1994), 113--170.
%
\bibitem{Fel99}
A.~A.~Felikson,
{\it Coxeter decompositions of spherical tetrahedra},
preprint, 99--053, Bielefeld, 1999.
%
\bibitem{GGM01}
F.~W.~Gehring, J.~P.~Gilman, and G.~J.~Martin,
{\it Kleinian groups with real parameters},
Commun. Contemp. Math. {\bf 3}, no.~2 (2001), 163--186.
%
\bibitem{GM94-1}
F.~W.~Gehring and  G.~J.~Martin,  {\it Chebyshev polynomials and
discrete groups}, Proc. of the Conf. on Complex Analysis (Tianjin,
1992), 114--125, Conf. Proc. Lecture Notes Anal., I, Internat.
Press, Cambridge, MA, 1994.
%
\bibitem{GM94-2}
F.~W.~Gehring and  G.~J.~Martin,
{\it On the minimal volume hyperbolic $3$-orbifold},
Math. Res. Lett. {\bf 1} (1994), no. 1, 107--114.
%
\bibitem{Kli89}
E.~Klimenko,  {\it Discrete  groups   in   three-dimensional
Lo\-ba\-chevsky space generated  by  two  rotations},  Siberian
Math.~J. {\bf 30} (1989), no.~1, 95--100.
%
\bibitem{Kli90}
E.~Klimenko, {\it Some remarks on subgroups of\/ ${\rm PSL}(2,{\mathbb C})$},
Q \& A
in General Topology {\bf 8} (1990), no.~2, 371--381.
%
\bibitem{Kli01}
E.~Klimenko,  {\it Some examples of discrete groups and
hyperbolic orbifolds of infinite volume},
J. Lie Theory {\bf 11} (2001), 491--503.
%
\bibitem{KK02}
E.~Klimenko and N.~Kopteva, {\it Discreteness criteria for
$\cal RP$ groups},
Israel J. Math. {\bf 128} (2002), 247--265.
%
\bibitem{KK04}
E.~Klimenko and N.~Kopteva, {\it All discrete $\mathcal{RP}$ groups
whose generators have real traces},
Int. J. Algebra Comput. {\bf 15} (2005), no.~3, 577--618.
%
\bibitem{KK04-1}
E.~Klimenko and N.~Kopteva, {\it Discrete ${\cal RP}$ groups with
a parabolic generator},
Sib. Math. J. {\bf 46} (2005), no.~6, 1069--1076.
%
\bibitem{KK05}
E.~Klimenko and N.~Kopteva, {\it Kleinian orbifolds
uniformized by ${\cal RP}$ groups with
an elliptic and a hyperbolic generators}, 2005, preprint.
%
\bibitem{MM03}
C.~Maclachlan and G.~J.~Martin,
{\it All Kleinian groups with two elliptic generators whose
commutator is elliptic},
Math. Proc. Camb. Philos. Soc. {\bf 135} (2003), no.3, 413-420.
%
\bibitem{Mat82}
J.~P.~Matelski, {\it The classification of discrete 2-generator
subgroups of ${\rm PSL}(2,{\mathbb R})$}, Israel J.~Math. {\bf 42} (1982),
no.~4, 309--317.
%
\bibitem{Vin85}
E.~B.~Vinberg,
{\it Hyperbolic reflection groups},
Russian Math. Surveys {\bf 40} (1985), 31--75.
\end{thebibliography}
\end{document}